\pdfoutput=1

\DeclareSymbolFont{AMSb}{U}{msb}{m}{n}
\documentclass[11pt,a4paper,leqno,noamsfonts]{amsart}

   \makeatletter
   \renewcommand\@biblabel[1]{#1.}
   \makeatother
   \usepackage{multirow}

\usepackage{graphicx}
\usepackage{subcaption}
\usepackage{pythontex}

\usepackage{listings} 

\usepackage[bookmarks=false]{hyperref}

\newtheorem{manualconjectureinner}{Conjecture}
\newenvironment{manualconjecture}[1]{%
  \IfBlankTF{#1}
    {\renewcommand{\themanualconjectureinner}{\unskip}}
    {\renewcommand\themanualconjectureinner{#1}}%
  \manualconjectureinner
}{\endmanualconjectureinner}

\usepackage{listings}
\usepackage{color}

\definecolor{dkgreen}{rgb}{0,0.6,0}
\definecolor{gray}{rgb}{0.5,0.5,0.5}
\definecolor{mauve}{rgb}{0.58,0,0.82}

\usepackage{mathtools}

\usepackage[english]{babel}
\usepackage[dvipsnames]{xcolor}
\usepackage{graphicx,pifont,soul,physics} 
\usepackage[utopia]{mathdesign}
\usepackage[a4paper,top=3cm,bottom=3cm,
           left=3cm,right=3cm,marginparwidth=60pt]{geometry}

\usepackage{fancyhdr}
\usepackage{float}

\usepackage[utf8]{inputenc}
\usepackage{braket,caption,comment,mathtools,stmaryrd}
\usepackage[usestackEOL]{stackengine}
\usepackage{multirow,booktabs,microtype,relsize}
\usepackage[bookmarks=false]{hyperref}%
      \hypersetup{colorlinks,%
           citecolor=britishracinggreen,%
            filecolor=black,%
            linkcolor=cobalt,%
            urlcolor=cornellred}
      \numberwithin{equation}{section}
\usepackage[capitalise]{cleveref}
\usepackage[colorinlistoftodos]{todonotes}

    \usepackage[bookmarks=false]{hyperref}
\definecolor{antiquewhite}{rgb}{0.98, 0.92, 0.84}
\definecolor{buff}{rgb}{0.94, 0.86, 0.51}
\definecolor{palecopper}{rgb}{0.85, 0.54, 0.4}
\definecolor{fluorescentyellow}{rgb}{0.8, 1.0, 0.0}
\definecolor{bole}{rgb}{0.47, 0.27, 0.23}

\usepackage{amsmath,float}
\usetikzlibrary{decorations.markings}

\usepackage[vcentermath]{youngtab}
\input xy
\xyoption{all}

\usepackage{pgfkeys}
\usepackage{pgfopts}
\usepackage{ytableau}

\definecolor{cornellred}{rgb}{0.7, 0.11, 0.11}
\definecolor{britishracinggreen}{rgb}{0.0, 0.26, 0.15}
\definecolor{cobalt}{rgb}{0.0, 0.28, 0.67}
\DeclareSymbolFont{usualmathcal}{OMS}{cmsy}{m}{n}
\DeclareSymbolFontAlphabet{\mathcal}{usualmathcal}

\newcommand{\BA}{{\mathbb{A}}}

\newcommand{\BC}{{\mathbb{C}}}

\DeclareMathOperator{\Hilb}{Hilb}

\DeclareMathOperator{\colength}{colength}

\DeclareFontFamily{OT1}{rsfs}{}
\DeclareFontShape{OT1}{rsfs}{n}{it}{<-> rsfs10}{}
\DeclareMathAlphabet{\curly}{OT1}{rsfs}{n}{it}
\renewcommand\hom{\mathscr{H}\kern-0.3em\mathit{om}}

\newcommand\Hom{\operatorname{Hom}}

\DeclareMathOperator{\lHom}{\mathscr{H}\kern-0.3em\mathit{om}}
\DeclareMathOperator{\RRlHom}{\mathbf{R}\kern-0.025em\mathscr{H}\kern-0.3em\mathit{om}}

\DeclareMathOperator{\lExt}{{\mathscr{E}\kern-0.2em\mathit{xt}}}



\usepackage[all]{xy}
\usepackage{tikz}
\usepackage{tikz-cd}
\usepackage{adjustbox}
\usepackage{rotating}
\usepackage{comment}

\usetikzlibrary{matrix,shapes,intersections,arrows,decorations.pathmorphing}
\tikzset{commutative diagrams/arrow style=math font}
\tikzset{commutative diagrams/.cd,
mysymbol/.style={start anchor=center,end anchor=center,draw=none}}

\tikzset{
shift up/.style={
to path={([yshift=#1]\tikztostart.east) -- ([yshift=#1]\tikztotarget.west) \tikztonodes}
}
}

\theoremstyle{definition}

\newtheorem*{lemma*}{Lemma}
\newtheorem*{theorem*}{Theorem}
\newtheorem*{example*}{Example}
\newtheorem*{fact*}{Fact}
\newtheorem*{notation*}{Notation}
\newtheorem*{definition*}{Definition}
\newtheorem*{prop*}{Proposition}
\newtheorem*{remark*}{Remark}
\newtheorem*{corollary*}{Corollary}

\newtheorem*{conventions*}{Conventions}

\newtheorem{definition}{Definition}[section]

\newtheorem{example}[definition]{Example}

\newtheoremstyle{thm} 
        {3mm}
        {3mm}
        {\slshape}
        {0mm}
        {\bfseries}
        {.}
        {1mm}
        {}
\theoremstyle{thm}

\newtheorem{conj}{Conjecture}

\newtheoremstyle{ex} 
        {3mm}
        {3mm}
        {}
        {0mm}
        {\scshape}
        {.}
        {1mm}
        {}
\theoremstyle{ex}

\newtheoremstyle{sol} 
        {3mm}
        {3mm}
        {}
        {0mm}
        {\scshape}
        {.}
        {1mm}
        {}
\theoremstyle{sol}

\usepackage{tikz}
\usepackage{xparse}

\newcount\tableauRow
\newcount\tableauCol

\newenvironment{Tableau}[1]{%
  \tikzpicture[scale=0.5,draw/.append style={thick,black},
                      baseline=(current bounding box.center)]
    \tableauRow=-1.5
    \foreach \Row in {#1} {
       \tableauCol=0.5
       \foreach\k in \Row {
         \draw[thin](\the\tableauCol,\the\tableauRow)rectangle++(1,1);
         \draw[thin](\the\tableauCol,\the\tableauRow)+(0.5,0.5)node{$\k$};

                  \global\advance\tableauCol by 1
       }
       \global\advance\tableauRow by -1
    }
}{\endtikzpicture}

\newtheorem*{Acknowledgments*}{Acknowledgments}
\DeclareMathAlphabet\BCal{OMS}{cmsy}{b}{n}

\address{Centre for Mathematical Sciences, University of Cambridge, Wilberforce Road, CB3 0WA, Cambridge, United Kingdom}

\title[A monotonicity conjecture]{A monotonicity conjecture for the local maximal singularity of the Hilbert scheme of points}

\author{Alexia Ascott}
\email{aga34@cam.ac.uk (Alexia Ascott)}
\author{Fatemeh Rezaee}
\email{fr414@cam.ac.uk (Fatemeh Rezaee)}

\address{Department of Mathematics, Stockholm University, Stockholm, Sweden}
\author{Zhichen Zhou}
\email{zhichen.zhou@math.su.se (Zhichen Zhou)}

\begin{document}
\maketitle
\begin{abstract} 
The Brian\c{c}on-Iarrobino conjecture predicts the maximum singularity of the Hilbert scheme of a tetrahedral number of points. As for the maximal singularities of the Hilbert scheme of a non-tetrahedral number of points, the second named author gave some separate conjectural necessary and sufficient conditions. In this paper, we provide a conjectural sufficient condition for the necessary condition, and propose a monotonicity conjecture which predicts that for a fixed colength $l$, the maximal dimension of the tangent space over all the Borel-fixed ideals of colength $l$ is increasing with respect to the smallest pure exponent of the ideal. 
\end{abstract}

{\hypersetup{linkcolor=black}
\tableofcontents}

\section{Introduction} 
Hilbert schemes of points are natural and crucial spaces in algebraic geometry. However, they are known to be highly singular. There has been tremendous growth in research to understand the singularities and behavior of these spaces.  For example, for Murphy's law for the Hilbert scheme, see \cite{Vakil06}, and for a proof of a variant for the Hilbert schemes of points,  see \cite{Jelisiejew20}; for the Hilbert scheme of points on threefolds, see \cite{Ramkumar-Sammartano,Jelisiejew-Ramkumar-Sammartano24, HuX, KatzS}; for the reducibility, see \cite{Iarrobino72, Hilb8, Jardim, Hilb_11}; for rational singularities, see \cite{Ramkumar-Sammartano24}; for irrational components of the Hilbert scheme of points, see \cite{Farkas-Pandharipande-Sammartano24}; for the parity conjecture, see \cite{PandharipandeSlides,Ramkumar-Sammartano1,GGGL23}; for other works related to enumerative geometry, see \cite{BBS, BFHilb,JKSCounterBehrend, MNOP1, Nesterov25, PT, RicolfiSign}.

Because of the extremely singular nature of the Hilbert schemes, measuring the singularity is of great importance.  Let $N\geq 3$. The dimension of the tangent space to $\mathrm{Hilb}^l(\BA^N)$ at each point measures how singular the point is. If $l$ is a tetrahedral number ${N-1+k\choose N}$,  Brian\c{c}on and Iarrobino predicted the most singular points of the Hilbert scheme, which was recently resolved in three dimensions in \cite{Mackenzie-Rezaee1}:

\begin{conj}[\cite{Bri-Iar}] \label{1978Conj} The most singular point of $\Hilb^{{N-1+k \choose N}}(\BA^N)$ corresponds to the ideal
 $\mathfrak{m}^k$, where $\mathfrak{m}$ is the maximal ideal in $\BC[x_1,x_2,\dots,x_N]$.
\end{conj}

 To generalize this  to non-tetrahedral numbers, investigating the relationship between the shape of the ideals corresponding to the points of the Hilbert scheme and their level of singularity was suggested in \cite{Rezaee-23-Conjectures}. In this article, we further generalize the necessary condition in \cite{Rezaee-23-Conjectures} using a novel pattern in the shape of the \textit{locally-maximal singularities} (see Definition \ref{def: local-max}).
 
First, we recall the necessary condition.

\begin{manualconjecture}{B}[Necessary condition for maximal singularity]  \cite[Conjecture B]{Rezaee-23-Conjectures}
\label{necessaryCondition} Let $A=\BC[x_1,x_2,\ldots,x_N]$. For $N\geq 3$,   let $I$ be a 0-dimensional Borel-fixed ideal of colength $l$ in $A$  given by
   \begin{align*}
       I=(x_1^{m_1},x_2^{m_2}, \ldots,x_N^{m_N},\text{mixed monomial generators}),
   \end{align*}
   where  $m_1\leq m_2\leq\ldots\leq m_N$. 
   If ${N-1+k \choose N}\leq l<{N+k\choose N}$, and $I$ correspond to the most singular point of $\mathrm{Hilb}^{l}(\BA^N)$, then $m_1=k$.
\end{manualconjecture}

We propose a stronger version of this conjecture. Before stating the conjecture, we introduce some notations and definitions.

For any Borel-fixed ideal $I=(x_1^{m_1},x_2^{m_2},\ldots,x_N^{m_N},\text{mixed terms})$, we have $m_1\leq m_2 \leq \ldots\leq m_N$. We recall that Borel-fixed ideals are among the candidates for obtaining the maximal dimension of the tangent space.

\begin{definition}[locally-maximal singularity or $m_1$-maximal singularity]\label{def: local-max}  Fix $l={k+N-1 \choose N}+\Delta$ and $m_1$.  We say that a Borel-fixed ideal $I$ of colength $l$ has a \textit{locally-maximal singularity} or \textit{$m_1$-maximal singularity} if it has the maximal dimension tangent space among all the Borel-fixed ideals of colength $l={k+N-1 \choose N}+\Delta$ and the exponent of $x_1$ being $m_1$. As mentioned above, we denote such a maximal dimension by $T_{\max,m_1}(l)$.
\end{definition}

Now, we have the following main statements.

\begin{manualconjecture}{C}[Monotonicity conjecture]\label{ConjC}
     For fixed $l$, the function $T_{\mathrm{max},m_1}(l)$ is increasing  in $m_1$. In particular, Conjecture \ref{necessaryCondition} holds.
\end{manualconjecture}

\begin{figure}[H]

  \subcaptionbox*{}[.98\linewidth]{
    \includegraphics[width=\linewidth]{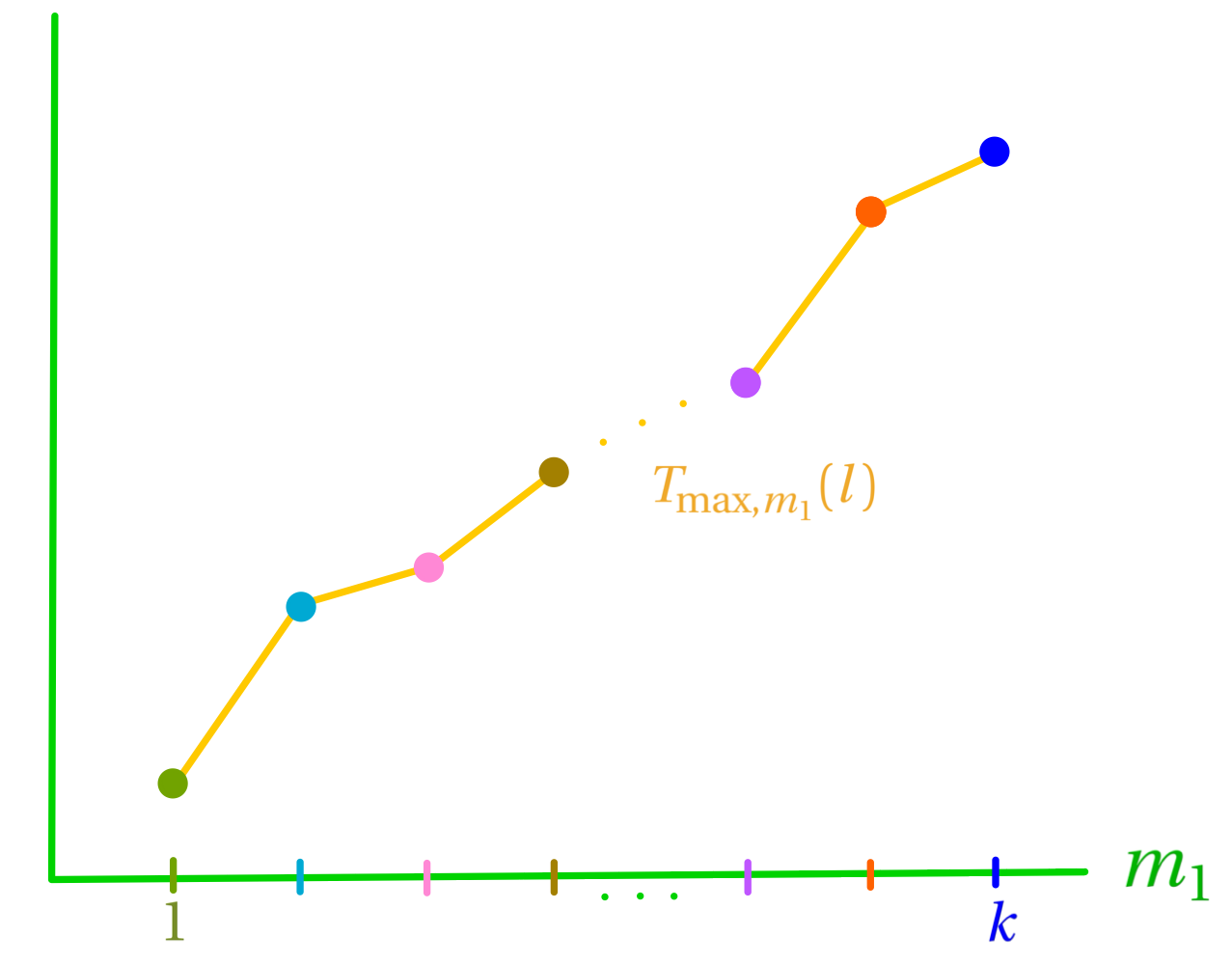}
  }
  \caption{For fixed $l$, $T_{\mathrm{max},m_1}(l)$ is increasing in $m_1$.}
  \label{Fig: Diagram}

\end{figure}

See Figure \ref{Fig: Diagram}.

\subsection*{Acknowledgements} This work is an outcome of the Cambridge Summer Research In Mathematics (SRIM) program 2024, where AA and ZZ were supervised by FR. The authors were supported by the ERC Advanced Grant MSAG. AA was also partially supported by St. John's College (Cambridge) funding, and FR was also supported by a UKRI grant EP/X032779/1. The authors thank Mark Gross for supporting the project, and the SRIM for making this work possible.  The authors also thank Tony Iarrobino  and Owen Mackenzie for helpful comments on the first version. We also thank Oscar Selby for his help with the code in Appendix \ref{App: proofApproach}.  
We acknowledge the use of  Macaulay2 \cite{M2} and Python.

\subsection*{Notation} We have the following maicn notations.
\begin{center}
    
         \begin{tabular}{ r l }

          $T(I)$:& The dimension of the  tangent space to the Hilbert scheme at the ideal $I$,\\& which is defined by $\mathrm{hom}(I,R/{I})$.\\

 $T_{\max,m_1}(l)$:& The maximal dimension of the tangent space for fixed $l$ and $m_1$\\& (see Definition \ref{def: local-max}).\\
   
 \end{tabular}
     \end{center}

\section{The monotonicity conjecture in examples}\label{Sec: ConjStatement}
In this brief section, we verify Conjecture \ref{ConjC} in some examples in the following table. In Appendix \ref{App: M2Codes}, we provide the codes we used to compute examples in this section. 
\subsection{Table} As evidence, in the following table, we confirm Conjecture \ref{ConjC}
for $N=3$ and $10\leq l\leq35$. Note that for $N=3$ and $m_1=1$, we have $T_{\max,m_1}(l)=3l$ (so this case is not included in the table).

\vspace{0.3cm}
\begin{tabular}{|p{0.86cm}|p{0.2cm}|p{0.35cm}|p{1.4cm}||||p{0.86cm}|p{0.2cm}|p{0.35cm}|p{1.4cm}||||p{0.86cm}|p{0.2cm}|p{0.35cm}|p{1.4cm}|} 
 \hline

 length ($l$)&k&$m_1$&$T_{\max,m_1}(l)$&length ($l$)&k&$m_1$&$T_{\max,m_1}(l)$&length ($l$)&k&$m_1$&$T_{\max,m_1}(l)$\\
 \hline\hline
 10&3&2&46& 21&4&2 &99&&&4& 190\\
  \cline{9-12}
&&3&60&&&3&129& 29&4&2 &141\\

 \cline{1-4}
 11&3&2&49&&&4& 153&&&3&183\\
 \cline{5-8}
&&3&63& 22&4&2 &104&&&4& 207\\
  \cline{9-12}
\cline{1-4}
 12&3&2&54&&&3&132& 30&4&2 &144\\ 
&&3&66&&&4& 156&&&3&186\\
 \cline{5-8} 
\cline{1-4}
 13&3&2&61& 23&4&2 &109&&&4& 210\\
  \cline{9-12}
&&3&69&&&3&141& 31&4&2 &149\\ 
\cline{1-4}
 14&3&2&64&&&4& 159&&&3&207\\
 \cline{5-8}
&&3&78& 24&4&2 &114&&&4& 217\\
  \cline{9-12} 
\cline{1-4}
15&3&2&69&&&3&148& 32&4&2 &156\\ 
&&3&81&&&4& 168&&&3&210\\
 \cline{5-8} 
\cline{1-4}
16&3&2&78& 25&4&2 &123&&&4& 240\\
  \cline{9-12} 
&&3&88&&&3&151& 33&4&2 &159\\ 
\cline{1-4}
17&3&2&81&&&4& 177&&&3&213\\
 \cline{5-8}  
&&3&99& 26&4&2 &126&&&4& 243\\
  \cline{9-12} 
\cline{1-4}
18&3&2&86&&&3&158& 34&4&2 &166\\
&&3&102&&&4& 180&&&3&216\\
 \cline{5-8} 
\cline{1-4}
19&3&2&89& 27&4&2 &129&&&4& 276\\
  \cline{9-12}  
&&3&123&&&3&161& 35&5&2 &169\\
\cline{1-4}
20&4&2&96&&&4& 187&&&3&225\\
 \cline{5-8}
&&3&126& 28&4&2 &134&&&4& 279\\
&&4&150&&&3&176&&&5& 315\\
\cline{9-12}
\cline{5-8}
\cline{1-4}

\end{tabular}
\vspace{0.01cm}

\appendix
\section{Macaulay2 codes}\label{App: M2Codes}
The codes we use are as follows. One of the most effective ways is to fix the colenght and also the number of (minimal) generators (from $2k+2$ to ${k+2 \choose 2}$) and list all the Borel-fixed ideals and the dimension of the tangent space at the corresponding point of the Hilbert scheme. In the following codes, we fix $\colength=35$ and the number of generators, $18$. The third line lists all Borel-fixed ideals $I$ of colength $35$ with $18$ (minimal) generators along with the dimension of the tangent space to the Hilbert scheme at the point $[I]$. 
\vspace{0.2cm}

\begin{verbatim}
i1 : loadPackage "StronglyStableIdeals";
i2 : S = stronglyStableIdeals(35, QQ[x,y,z,w]);
i3 : for I in S list (if numgens(I)=!=18 then continue;
(I,degree Hom(I,(ring I)/I)))
\end{verbatim}

\section{A potential approach to prove the conjectures}\label{App: proofApproach}

Let $I$ be a Borel-fixed ideal of colength $l$ in $\BC[x,y,z]$. Also, let $G$ be the (minimal) number of generators of $I$. By definition of the tangent space as $\Hom(I,R/I)$, to explicitly count the tangent vectors, one approach is to eliminate the vectors in the $Gl$-space of the potential maps from the generators of $I$ to $R/I$, which are not in the tangent space. We call such vectors \textit{zero vectors}. Spotting patterns among the set of zero vectors for fixed $k$, $l={k+2\choose 3}$ and $m_1$ (the exponent of $x$, with the convention of being the smallest pure exponent, as in the first paragraph of Section \ref{Sec: ConjStatement}) will lead to proving Conjecture \ref{ConjC}. 

To visualize the zero vectors, we consider a non-tetrahedral but straightforward example to stress that the method can be applied to any Borel-fixed ideal.

\begin{example} Let $I=(x,y,z^2)^2$. Using the codes, we can list the zero vectors illustrated in Figure \ref{Fig: zeroVectors}. Recall that this ideal has colenght $8$ and it has the maximal dimension tangent space, as considered in \cite{Sturmfels}. Using the zero vectors approach, we can calculate the dimension of the tangent space, $6\times 8-12=36$.

       \begin{figure}[h]
 \subcaptionbox*{(i)}[.35\linewidth]{
    \includegraphics[width=\linewidth]{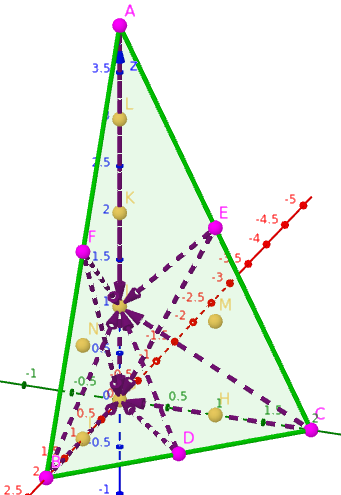}
  }
  \hskip15ex
  \subcaptionbox*{(ii)}[.28\linewidth]{
    \includegraphics[width=\linewidth]{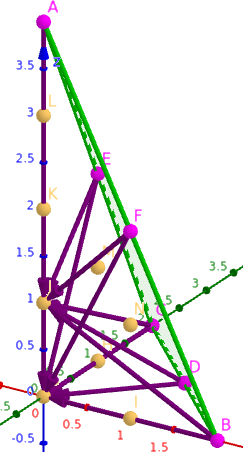}
  }
  \caption{12 zero vectors for $I=(x,y,z^2)^2$.}
  \label{Fig: zeroVectors}

\end{figure}
\end{example}

Here, we provide some codes (in Python) that help list the zero vectors of any Borel-fixed ideal in three variables. Using the approach in \cite{Ramkumar-Sammartano}, one can recognize the tangent vectors in 3D and then look at the complements as zero vectors. The following example calculates the dimension of the graded component of $T(I)$ of degree $\alpha \in \mathbb{Z}^3$ when $I=(x^2,y^3,z^3,xy,xz,yz^2,y^2z)$ and $\alpha = (0,2,-3)$.

\begin{verbatim}
import numpy as np
import sys
from scipy.ndimage import label
# np.set_printoptions(threshold=sys.maxsize)

# THESE ARE THE INPUT POINTS
input_points = 
np.array([(2,0,0),(0,3,0),(0,0,3),(1,1,0),(1,0,1),(0,1,2),(0,2,1)])
x_max, y_max, z_max = (np.max(input_points[:, 0]),
                       np.max(input_points[:, 1]),
                       np.max(input_points[:, 2]))

size = ((x_max + y_max + z_max + 3)*3)
output_matrix = np.zeros((size*2, size*2, size*2), dtype=int)
for point in input_points:
    output_matrix[point[0]+size:, point[1]+size:, point[2]+size:] = 1

output_matrix *= 2

# print(output_matrix[:, :, :])

# THIS IS THE INPUT VECTOR
vector = (0, 2, -3)
rolled = np.roll(output_matrix, vector[0], axis=0)
rolled = np.roll(rolled, vector[1], axis=1)
rolled = np.roll(rolled, vector[2], axis=2)

if vector[0] < 0:
    rolled[vector[0]:, :, :] = 2

if vector[0] > 0:
    rolled[:vector[0], :, :] = 2

if vector[1] < 0:
    rolled[:, vector[1]:, :] = 2

if vector[1] > 0:
    rolled[:, :vector[1], :] = 2

if vector[2] < 0:
    rolled[:, :, vector[2]:] = 2

if vector[2] > 0:
    rolled[:, :, :vector[2]] = 2

rolled //= 2

# print("This is the shifted matrix: \n", rolled)
combined = rolled + output_matrix
combined[np.where(combined > 1)] = 0
# print("This is the combined matrix: \n", combined)

structure = np.array([[[0, 0, 0],
                       [0, 1, 0],
                       [0, 0, 0]],
                      [[0, 1, 0],
                       [1, 1, 1],
                       [0, 1, 0]],
                      [[0, 0, 0],
                       [0, 1, 0],
                       [0, 0, 0]]
                      ])
labeled, ncomponents = label(combined, structure)
indices = np.indices(combined.shape).T[:, :, :, [1, 0]]
total = 0
for n in range(1, ncomponents + 1):
    if len(indices[labeled == n]) <= (x_max + y_max + z_max + 3) ** 2:
        total += 1

print(total)
\end{verbatim}

\bibliographystyle{amsplain-nodash}

\bibliography{bib}

\vspace{0.25cm}

\end{document}